\documentclass{amsart}

\usepackage{amsmath}
\usepackage{amsfonts}
\usepackage{amssymb}

\newcommand{\half}{\tfrac{1}{2}}
\newcommand{\summ}{\mathop{{\sum}^{\star}}}
\newcommand{\sumh}{\mathop{{\sum}^{p}}}

\numberwithin{equation}{section}

\newtheorem{theorem}{Theorem}[section]

\newtheorem{lemma}[theorem]{Lemma}
\newtheorem{corollary}[theorem]{Corollary}

\begin{document}

\title{Simultaneous non-vanishing of $GL(3)\times GL(2)$ and $GL(2)$ $L$-functions.}
\author{Rizwanur Khan}
\address{Mathematisches Institut, Georg-August Universit\"{a}t G\"{o}ttingen, Bunsenstra$\ss$e 3-5,
         D-37073 G\"{o}ttingen, Germany}
\email{rrkhan@uni-math.gwdg.de}

\begin{abstract}
Fix $g$ a Hecke-Maass form for $SL_3(\mathbb{Z})$. In the family of holomorphic newforms $f$ of fixed weight and large prime level $q$, we find the average value of the product $L(\half,g\times f)L(\half,f)$. From this we derive a result on the simultaneous non-vanishing of these $L$-functions at the central point.
\end{abstract}

\maketitle

\section{Introduction}

In the analytic theory of $L$-functions, extensive progress has been made for $GL(1)$ and $GL(2)$ $L$-functions towards understanding their behaviour at the centre of the critical strip.  However, results on mean values, non-vanishing, size, etc. of {\it higher} rank $L$-functions are relatively few and very desirable. In this paper we study the mean value at the central point of a product of a degree 6 and a degree 2 $L$-function, and deduce a result on their simultaneous non-vanishing.

Let $H_k^{\star}(N)$ denote the set of holomorphic cusp forms of  even weight $k$ and trivial nebentypus which are newforms of level $N$ in the sense of Atkin-Lehner theory \cite{atk}. Let $L(s,f)$ be the $L$-function attached to $f\in H_k^{\star}(N)$, normalized so that its functional equation relates values at $s$ and $1-s$. Kowalski, Michel and VanderKam \cite{kmv}, building on the work of Duke, Friedlander, and Iwaniec \cite{dfi2, dfi}, found amongst other things an asymptotic for the following fourth power mean value, as $q\to \infty$ amongst the primes:
\begin{align}
\label{fourthm} \sumh_{f\in H_2^{\star}(q)} L(\half,f)^4,
\end{align}
where $\sumh$ is the Petersson average defined in (\ref{harmonic}). In subsequent work, Kowalski, Michel and VanderKam \cite{kmv2, mic} generalized the fourth moment above by studying the second moment of $GL(2)\times GL(2)$ Rankin-Selberg $L$-functions. Fix $h$ a holomorphic Hecke cusp form for $SL_2(\mathbb{Z})$. One of their many results was an asymptotic for the following square mean value, as $q\to \infty$ amongst the primes:
\begin{align}
\sumh_{f\in H_k^{\star}(q)} |L(\half, h \times f)|^2,
\end{align}
for $k<12$. In this paper, we generalize (\ref{fourthm}) in a different direction. Fix $g$ a Hecke-Maass form for $SL_3(\mathbb{Z})$ which is unramified at infinity and let $\tilde{g}$ denote its dual.  Let  $L(s, g\times f)$ denote the $GL(3)\times GL(2)$ Rankin-Selberg convolution of $g$ with $f\in H_k^{\star}(N)$. This is defined in section \ref{autback} so that the central point equals $s=\half$. We prove
\begin{theorem}\label{main}
Follow the notation in section \ref{autback}. Let $\epsilon>0$. There exists $k_\epsilon>0$ such that for even $k>k_\epsilon$ and $q$ prime, we have
\begin{multline}
\label{mainline} \frac{12}{q(k-1)}\sumh_{f\in H_k^{\star}(q)} L(\half, g\times f) L(\half, f) \\
= \frac{L(1,g)L(1,\tilde{g})}{\zeta(2)}\Big(1+\frac{\widetilde{G}_1(\half)}{G_1(\half)}\Big) + O(q^{\frac{2\theta_1-1}{2}+\epsilon}+q^{\frac{2\theta_2-1}{4}+\epsilon}+q^{-\frac{1}{8}+\epsilon}),
\end{multline}
where $\theta_1$ and $\theta_2$ are bounds towards the Ramanujan conjecture for $GL(3)$ in the non-archimedian and archimedian aspects respectively, and the implied constant depends on $g$, $\epsilon$ and $k$.
\end{theorem}
\noindent For the sake of clarity, we do not specify $k_\epsilon$, but this can be done. The theorem gives an asymptotic for any non-trivial bounds towards the Ramanujan conjecture for $GL(3)$. If $g$ is self-dual then we may take $\theta_1 = 7/32$ and $\theta_2= 0$. Otherwise, the best bounds currently known are $\theta_1=\theta_2=5/14$. These facts are discussed further in section \ref{autback}. If $g$ is a self-dual form, it is known by the work of Soudry \cite{soud} to be the symmetric square lift of a Hecke-Maass form for $SL_2(\mathbb{Z})$. In this case $L(s, g\times f) L(s, f)$ is essentially a triple product $L$-function. 

The $GL(3)\times GL(2)$ $L$-functions are important and have been studied by other authors. In a recent breakthrough paper, Li \cite{li1} proved a subconvexity bound for $L(1/2,g\times f)$ when $g$ is self-dual and $f\in H^\star_k(1)$, in the $k$-aspect. She also considered twists by Hecke-Maass forms for $SL_2(\mathbb{Z})$ and proved subconvexity in the eigenvalue aspect. Blomer \cite{blom} proved a subconvexity bound for $L(1/2, g\times f)$ for $g$ self-dual and special Hecke-Maass forms $f$ for $\Gamma_0(q^2)$, where $q$ is prime, in the $q$-aspect. 

Questions on the simultaneous non-vanishing of two $L$-functions at the central point are also of interest. For example, in the paper \cite{kmv} mentioned above, the authors built on (\ref{fourthm}) to prove a result on the simultaneous non-vanishing of $L(1/2,f)$ and $L(1/2,f\times \chi)$, where $\chi$ is a fixed, non-quadratic primitive Dirichlet character. Simultaneous non-vanishing in the case when $\chi$ is quadratic was considered by Iwaniec and Sarnak \cite{iwasar}, in their work on Landau-Siegel zeros. Before her work on subconvexity, Li \cite{li2} studied the simultaneous non-vanishing of $GL(3)\times GL(2)$ and $GL(2)$ $L$-functions, in the $GL(2)$ family of Hecke-Maass forms for $SL_2(\mathbb{Z})$. We address this problem in the level aspect.

\begin{corollary}
For all prime $q$ and even $k$ larger than some constant depending on $g$, there exists $f\in H_k^{\star}(q)$ such that  $L(\half, g\times f) \neq 0$ and $ L(\half, f) \neq 0$. If $g$ is self-dual then $k$ must be larger than an absolute constant.
\end{corollary}
\proof
We need to show that the main term of Theorem \ref{main} is non-zero. If $g$ is self-dual then $\frac{\widetilde{G}_1(1/2)}{G_1(1/2)}=1$. Otherwise, by Stirling's asymptotic for the gamma function and the identity $\alpha_1+\alpha_2+\alpha_3=0$, we have
\begin{align}
\frac{\widetilde{G}_1(\half)}{G_1(\half)} = k^{-\alpha_1-\alpha_2-\alpha_3}(1+O_g(k^{-1})) = 1+O_g(k^{-1}),
\end{align}
which is non-zero for $k$ larger than a constant depending on $g$. We also have by \cite{jacsha} that $L(1,g)$ and $L(1,\tilde{g})$ are non-zero.
\endproof

{\it Notation.} Throughout, $\epsilon$ will denote an arbitrarily small positive constant, but not necessarily the same one from one occurrence to the next. Any implied constant may depend implicitly on $g$, $\epsilon$ and $k$. For real numbers $a,b>0$, we write $a \asymp b$ to mean $q^{-\epsilon}<a/b< q^{\epsilon}$.

\subsection{Sketch}\label{sket}

Here we give an imprecise outline of a part of the proof of Theorem \ref{main}. By the method of approximate function equations, we express $L(\half, g\times f)$ as a weighted Dirichlet sum of length about $q^{3/2}$ and $L(\half, f)$ as a weighted Dirichlet sum of length about $q^{1/2}$. A part of what we must calculate is similar to
\begin{align}
\label{out1} \frac{12}{q(k-1)}\sumh_{f\in H_k^{\star}(q)} \sum_{\substack{r^2n< q^{3/2}\\ m< q^{1/2}}} \frac{A(r,n)a_f(n) a_f(m)}{r\sqrt{nm}},
\end{align}
where $\sumh$ denotes a weighted sum as on the left hand side of (\ref{mainline}) and $A(r,n)$ and $a_f(n)$ are the Fourier coefficients of $g$ and $f$ respectively, suitably normalized. Precise definitions are made in the following sections. Applying a trace formula, (\ref{out1}) essentially equals
\begin{align}
\sum_{r,n< q^{1/2}} \frac{A(r,n)}{rn} +  \sum_{\substack{r^2n< q^{3/2}\\ m< q^{1/2}\\ c\ge 1}} \frac{A(r,n)}{r\sqrt{nm}} \frac{S(n,m,cq)}{cq}J_{k-1}\Big(\frac{4\pi\sqrt{nm}}{cq}\Big).
\end{align}
The first sum above essentially gives rise to the main term. We must show that the second sum falls into the error term. By the rapid decay of the $J$-Bessel function for small argument, when $k>k_{\epsilon}$, we see that we may assume that $c,r<q^{\epsilon}$. Let us assume $c=r=1$. Opening the Kloosterman sum, a part of what we must bound is
\begin{align}
\sum_{q^{1/2}< m< 2q^{1/2}} \summ_{h \bmod q} e(mh/q) \sum_{q^{3/2}< n< 2q^{3/2}} A(1,n) e(n\overline{h}/q).
\end{align}
(The $J$-Bessel function is roughly constant in this range of $n$ and $m$.) We apply the $GL(3)$ Voronoi summation formula to exchange the $n$-sum for another sum of length about $q^3/q^{3/2}=q^{3/2}$. A part of what we must bound is then
\begin{align}
 \sum_{q^{1/2}< m< 2q^{1/2}} \summ_{h \bmod q} e(mh/q) \sum_{q^{3/2}< n< 2q^{3/2}} A(n,1) S(n,h,q)
\end{align}
We have $\summ_{h \bmod q} e(mh/q) S(n,h,q) \approx q e(n\overline{m}/q)$. By reciprocity (the Chinese Remainder Theorem), we have $e(n\overline{m}/q)=e(n/mq)e(-n\overline{q}/m)\approx e(-n\overline{q}/m)$, since $n\approx q^{3/2} \approx mq$. Thus we must bound
\begin{align}
\sum_{q^{1/2}< m< 2q^{1/2}}  \sum_{q^{3/2}< n< 2q^{3/2}} A(n,1) e(-n\overline{q}/m).
\end{align}
The new modulus $m$ of the exponential is much smaller than the original modulus $q$. We apply the $GL(3)$ Voronoi summation formula again, to exchange the $n$-sum for another sum of length about $m^3/q^{3/2}\approx 1$. We must bound
\begin{align}
\sum_{q^{1/2}< m< 2q^{1/2}} \sum_{ n< q^{\epsilon}} A(1,n) S(-n,q,m).
\end{align}
We bound this sum absolutely, using Weil's bound for the Kloosterman sum.

\subsection{Automorphic forms and $L$-functions}\label{autback}

Every $f\in H_k^{\star}(q)$ has a Fourier expansion of the type
\begin{align}
f(z) = \sum_{n=1}^{\infty} a_f(n) n^{\frac{k-1}{2}} e(nz)
\end{align}
for $\Im z>0$, where $e(z)=e^{2\pi i z}$, $a_f(n)\in \mathbb{R}$ and $a_f(1)=1$. The coefficients $a_f(n)$ satisfy the multiplicative relation
\begin{align}
\label{hmult} a_f(n)a_f(m) = \sum_{\substack{d|(n,m)\\(d,q)=1}} a_f\Big(\frac{nm}{d^2}\Big)
\end{align}
and Deligne's bound $a_f(n)\le d(n)$. The $L$-function associated to $f$ is defined as
\begin{align}
L(s,f) = \sum_{n=1}^{\infty} \frac{a_f(n)}{n^s}
\end{align}
for $\Re(s)>1$. This satisfies the functional equation
\begin{align}
\label{funct1} q^{\frac{s}{2}}   G_2(s) L(s,f) = \epsilon_f q^{\frac{1-s}{2}}  G_2(1-s) L(1-s,f),
\end{align}
where 
\begin{align}
G_2(s) = \pi^{-s}  \Gamma\Big( \frac{s+ \frac{k-1}{2}}{2} \Big)\Gamma\Big( \frac{s+ \frac{k+1}{2}}{2} \Big)
\end{align}
and $\epsilon_f = - i^k  a_f(q) q^{\frac{1}{2}}= \pm 1$. The left hand side of (\ref{funct1}) analytically continues to an entire function. The facts above can be found in \cite{iwaniec}.

We fix a Hecke-Maass form of type $(\nu_1,\nu_2)$ for $SL_3(\mathbb{Z})$. We refer to \cite{gold}, especially Chapter 6, and follow its notation. We write $A(n,m)$ for the Fourier coefficients of $g$ in the Fourier expansion (6.2.1) of \cite{gold}, normalized so that $A(1,1)=1$. The $L$-functions associated to $g$ and its dual $\tilde{g}$ are defined as
\begin{align}
L(s,g) = \sum_{n=1}^{\infty} \frac{A(n,1)}{n^s}, \hspace{0.2in} L(s,\tilde{g}) =   \sum_{n=1}^{\infty} \frac{\overline{A(n,1)}}{n^s}= \sum_{n=1}^{\infty} \frac{A(1,n)}{n^s}
\end{align}
for $\Re(s)>1$.
We  have the Hecke relation
\begin{align}
\label{hecke3} A(n,m) = \sum_{d|(n,m)} \mu(d) A\Big(\frac{n}{d},1\Big)A\Big(1,\frac{m}{d}\Big),
\end{align}
and if $(n_1m_1,n_2m_2)=1$, we have
\begin{align}
A(n_1n_2, m_1m_2)= A(n_1,m_1)A(n_2,m_2).
\end{align}
Let 
\begin{align}
\label{g1} \alpha_1= -1+\nu_1+2\nu_2, \ \alpha_2 = \nu_1-\nu_2, \ \alpha_3= -2\nu_1 -\nu_2 + 1.
\end{align}
Suppose that we have the following bounds towards the Ramanujan conjecture, at the finite and infinite places respectively:
\begin{align}
\label{rama} |A(n,1)| \ll n^{\theta_1+\epsilon}, \hspace{0.2in} |\Re(\alpha_1)|, |\Re(\alpha_2)|, |\Re(\alpha_3)| < \theta_2 +\epsilon
\end{align}
for some $0\le \theta_1,\theta_2\le \half$. While we expect $\theta_1=\theta_2=0$, the best bounds currently known are $\theta_1=\theta_2=\frac{5}{14}$. This is implicit in the work of Kim and Sarnak \cite{kimsar} and was noticed by Blomer and Brumley \cite{blobru}. In the case that $g$ is self-dual, we know more. If $g$ is a self-dual form, it is the symmetric square lift of a Hecke-Maass form for $SL_2(\mathbb{Z})$. By Kim and Sarnak's bounds \cite{kimsar}, we may take $\theta_1=\frac{7}{32}$ and since there are no exceptional eigenvalues for $SL_2(\mathbb{Z})$, we can take $\theta_2 =0$. By (\ref{rama}) and (\ref{hecke3}) we have
\begin{align}
A(n,m)\ll (nm)^{\theta_1+\epsilon}.
\end{align}
By (\ref{rama}) and Rankin-Selberg theory we have (cf. \cite{blom} for a proof):
\begin{align}
\sum_{n\le x} |A(na,b)|^2 \ll x(ab)^{2\theta_1+\epsilon}.
\end{align}
This together with (\ref{hecke3}) and the Cauchy-Schwarz inequality yields
\begin{align}
\label{avg} \sum_{\substack{n<x\\m<y}} |A(na,mb)| \ll (xy)^{1+\epsilon} (ab)^{\theta_1 +\epsilon}.
\end{align}

The Rankin-Selberg $L$-function $L(s, g\times f)$ is defined as
\begin{align}
L(s, g\times f)= \sum_{n,r\ge 1} \frac{A(r,n)a_f(n)}{(r^2n)^s}
\end{align}
for $\Re(s)>1$. It satisfies the functional equation
\begin{align}
\label{funct2}  q^{\frac{3s}{2}} G_1(s)  L(s, g\times f) = \epsilon_{g\times f} q^{\frac{3(1-s)}{2}}  \widetilde{G}_1(1-s)  L(1-s,  \tilde{g} \times f),
\end{align}
where 
\begin{align}
\epsilon_{g\times f}=(\epsilon_f)^3=\epsilon_f,
\end{align}
\begin{multline}
G_1(s) = \pi^{-3s} \Gamma\Big( \frac{s+\frac{k+1}{2} + \alpha_1}{2} \Big)\Gamma\Big( \frac{s+\frac{k+1}{2} + \alpha_2 }{2} \Big)\Gamma\Big( \frac{s+\frac{k+1}{2} + \alpha_3}{2} \Big)\\
 \times \Gamma\Big( \frac{s+ \frac{k-1}{2} + \alpha_1}{2} \Big)\Gamma\Big( \frac{s+ \frac{k- 1}{2} + \alpha_2}{2} \Big)\Gamma\Big( \frac{s+\frac{k-1}{2} + \alpha_3}{2} \Big),
\end{multline}
and
\begin{multline}
\widetilde{G}_1(s) = \pi^{-3s} \Gamma\Big( \frac{s+\frac{k+1}{2} - \alpha_1}{2} \Big)\Gamma\Big( \frac{s+\frac{k+1}{2} - \alpha_2 }{2} \Big)\Gamma\Big( \frac{s+\frac{k+1}{2} - \alpha_3}{2} \Big)\\
 \times \Gamma\Big( \frac{s+ \frac{k-1}{2} - \alpha_1}{2} \Big)\Gamma\Big( \frac{s+ \frac{k-1}{2} - \alpha_2}{2} \Big)\Gamma\Big( \frac{s+\frac{k-1}{2} - \alpha_3}{2} \Big).
\end{multline}
The left hand side of (\ref{funct2}) analytically continues to an entire function. To study these $L$-functions at $s=1/2$, we first express the central values as Dirichlet-type sums, using a standard technique.

\begin{lemma}{\bf Approximate functional equations} 

(i) Let 
\begin{align}
\label{v1def} V_1(x)= \frac{1}{2\pi i} \int_{(\sigma)} x^{-s} \frac{G_1(\half + s)}{G_1(\half)} \frac{ds}{s}, \hspace{0.3in} \widetilde{V}_1(x)= \frac{1}{2\pi i} \int_{(\sigma)} x^{-s} \frac{\widetilde{G}_1(\half + s)}{G_1(\half)} \frac{ds}{s},
\end{align}
for $x, \sigma>0$.
We have
\begin{align}
\label{afe1} L(\half,f\times g) = \sum_{n,r\ge 1} \frac{a_f(n)A(r,n)}{r\sqrt{n}}V_1\Big(\frac{r^2n}{q^{3/2}}\Big)-i^kq^{1/2}a_f(q)\sum_{n,r\ge 1} \frac{a_f(n)A(n,r)}{r\sqrt{n}}\widetilde{V}_1\Big(\frac{r^2n}{q^{3/2}}\Big).
\end{align}

(ii) Let 
\begin{align}
V_2(x)= \frac{1}{2\pi i} \int_{(\sigma)} x^{-s} \frac{G_2(\half + s)}{G_2(\half)} \frac{ds}{s}
\end{align}
for $x,\sigma >0$.
We have
\begin{align}
\label{afe2} L(\half,f) = \sum_{m\ge 1} \frac{a_f(m)}{\sqrt{m}}V_2\Big(\frac{m}{q^{1/2}}\Big)-i^k q^{1/2}a_f(q)\sum_{m\ge 1} \frac{a_f(m)}{\sqrt{m}}V_2\Big(\frac{m}{q^{1/2}}\Big).
\end{align}
For any $A>0$ and integer $B\ge 0$ we have that
\begin{align}
\label{vbound} V_1^{(B)}(x), \widetilde{V}_1^{(B)}(x), V_2^{(B)}(x) \ll_B x^B (1+x)^{-A},
\end{align}
so that the sums in (\ref{afe1}) and (\ref{afe2}) are essentially supported on  $r^2n<q^{3/2+\epsilon}$ and $m<q^{1/2+\epsilon}$.
\end{lemma}
\proof
The proof is a slight modification of the proof Theorem 5.3 of \cite{iwakow}, since the gamma factors $G_1(s)$ and $\widetilde{G}_1(s)$ may not be identical.
\endproof

We make the following simple observation for later use.
\begin{lemma}\label{vnote}
For $f_0 \in H_k^{\star}(1)$ and $\delta>0$ we have
\begin{align}
\label{vnote1} \sum_{n,r\ge 1} \frac{A(r,n)a_{f_0}(n)}{{r\sqrt{n}}} V_1\Big(\frac{r^2n}{q^{\delta}}\Big) \ll_{g\times f_0} 1
\end{align}
and
\begin{align}
\label{vnote2} \sum_{m\ge 1} \frac{a_{f_0}(m)}{{\sqrt{m}}} V_2(\frac{m}{q^{\delta}}) \ll_{f_0} 1.
\end{align}
\end{lemma}
\proof By definition (\ref{v1def}) we have that (\ref{vnote1}) equals
\begin{align}
\frac{1}{2\pi i } \int_{(2)} L(g\times f_0, \half +s) q^{\delta s} \frac{G_1(\half + s)}{G_1(\half)} \frac{ds}{s}.
\end{align}
We may move the line of integration to the left of the imaginary axis to get the required bound. One can simply take the convexity bound for the $L$-function since $g$ and $f_0$ are fixed. The proof of (\ref{vnote2}) is similar.
\endproof

\subsection{Trace formula}

As usual,
\begin{align}
S(n,m;c)= \summ_{h \bmod c} e\Big(\frac{nh+m\overline{h}}{c}\Big)
\end{align}
will denote the Kloosterman sum, where $\star$ indicates that the summation is restricted to $(h,c)=1$ and where $h\overline{h}\equiv 1 \bmod c$. We have Weil's estimate
\begin{align}
\label{weil} |S(n,m;c)| \le (n,m,c)^{\half} c^{\half} d(c).
\end{align}

$J_{k-1}(x)$ will denote the $J$-Bessel function. We will need the following estimates which can be found in \cite{gr} and \cite{wat}.
\begin{lemma}
Let $x>0$. We have
\begin{align}
\label{jbound} J_{k-1}(x) \ll \text{min}( x^{k-1}, x^{-1/2}).
\end{align}
We have
\begin{align}
\label{jlong} J_{k-1}(x) = \Re \Big( e^{ix}\mathcal{J}(x) \Big),
\end{align}
where $\mathcal{J}(x)$ is a smooth function depending on $k$, which satisfies
\begin{align}
\label{jbound2} \mathcal{J}^{(B)}(x)\ll_{B} x^{-B+1}(1+x)^{-3/2}
\end{align}
for any integer $B\ge 0$.
\end{lemma}

For any complex numbers $\lambda_f$, define the weighted sum
\begin{align}
\label{harmonic} \sumh_{f\in H_k^{\star}(N)}  \lambda_f = \sum_{f\in H_k^{\star}(N)} \frac{\lambda_f}{\zeta(2)^{-1} L(1, \text{sym}^2 f)},
\end{align}
where $ L(s, \text{sym}^2 f)$ denotes the $L$-function of the symmetric-square lift of $f$. The arithmetic weights above occur naturally in the Petersson trace formula (\ref{ptf}) and the following trace formula for newforms. Define
\begin{align}
\Delta^{\star}_{k,N}(n,m)=\frac{12}{N(k-1)}& \sumh_{f\in H_k^{\star}(N)} a_f(n)a_f(m).
\end{align}

\begin{lemma}{\bf Trace formula.} 

(i) We have
\begin{align}
\label{ptf} \Delta^{\star}_{k,1}(n,m)= \delta(n,m)  + 2\pi i^k \sum_{c\ge 1} \frac{S(n,m;c)}{c}J_{k-1}\Big(\frac{4\pi\sqrt{nm}}{c}\Big),
\end{align}
where $\delta(n,m)$ equals $1$ if $n=m$ and $0$ otherwise.

(ii) Let $q$ be a prime. If $(m,q)=1$ and $q^2\nmid n$ then
\begin{multline}
 \label{trace}\Delta^{\star}_{k,q}(n,m)  = \delta(n,m) + 2\pi i^k \sum_{c\ge 1} \frac{S(n,m;cq)}{cq}J_{k-1}\Big(\frac{4\pi\sqrt{nm}}{cq}\Big)\\
 -\frac{1}{q[\Gamma_0(1):\Gamma_0((n,q))]}\sum_{i=0}^{\infty} q^{-i} \Delta^{\star}_{k,1}(n, mq^{2i}).
 \end{multline}
\end{lemma}
\proof
See Proposition 2.8 of \cite{ils}. Note the different normalization there.
\endproof
\noindent The left hand side and the right hand side of the trace formula are sometimes called the spectral side and the arithmetic side respectively. We note that the last line of (\ref{trace}) is $\ll q^{-1+\epsilon}(nm)^{\epsilon}$. If $q|n$ then the last line is actually $\ll q^{-2+\epsilon}(nm)^{\epsilon}$ since $[\Gamma_0(1):\Gamma_0(q)]>q$.
We also note, using (\ref{jbound}), that for $nm<q^{2-\epsilon}$, we have
\begin{align}
\label{tracenote} \Delta^{\star}_{k,q}(n,m) = \delta(n,m) + O(q^{-100})
\end{align}
for $k$ large enough.

\subsection{Summation Formula}

The GL(3) Voronoi summation formula (\ref{voron}) was proved by Miller and Schmid \cite{milsch}. Goldfeld and Li \cite{golli} later gave another proof, and we follow their presentation. The asymptotic (\ref{further}) is due to Ivic \cite{ivi} and Li \cite{li1}, but we follow the presentation of \cite{blom}. 

\begin{lemma}{\bf GL(3) Voronoi Summation}
Let $\psi$ be a smooth, compactly supported function on the positive real numbers and $(b,d)=1$. We have
\begin{align}
\label{voron} \sum_{n\ge 1} A(r,n) e\Big(\frac{n\overline{b}}{d}\Big) \psi\left( \frac{n}{N} \right) =
 \sum_{\pm} \frac{\pi^{\frac{3}{2}}}{2} d\sum_{\substack{n\ge 1\\ l | dr}} \frac{A(n,l)}{nl}S\Big(rb,\pm n;\frac{dr}{l}\Big)\Psi^{\pm}\Big(\frac{n}{d^3N^{-1}rl^{-2}}\Big),
\end{align}
where we define
\begin{align}
&\Psi^{\pm}(X)=  X \frac{1}{2\pi i}\int_{(\sigma)} (\pi^3 X)^{-s} H^{\pm}(s) \widetilde{\psi}(1-s) ds ,\\
&\nonumber H^{\pm}(s)= \frac{\Gamma\Big(\frac{s+\alpha_1}{2}\Big)\Gamma\Big(\frac{s+\alpha_2}{2}\Big)\Gamma\Big(\frac{s+\alpha_3}{2}\Big)}{\Gamma\Big(\frac{1-s-\alpha_1}{2}\Big)\Gamma\Big(\frac{1-s-\alpha_2}{2}\Big)\Gamma\Big(\frac{1-s-\alpha_3}{2}\Big)} 
\mp i \frac{\Gamma\Big(\frac{1+s+\alpha_1}{2}\Big)\Gamma\Big(\frac{1+s+\alpha_2}{2}\Big)\Gamma\Big(\frac{1+s+\alpha_3}{2}\Big)}{\Gamma\Big(\frac{2-s-\alpha_1}{2}\Big)\Gamma\Big(\frac{2-s-\alpha_2}{2}\Big)\Gamma\Big(\frac{2-s-\alpha_3}{2}\Big)}
\end{align}
for $\sigma > \theta_2$, where $\alpha_i$ are defined as in (\ref{g1}) and $ \widetilde{\psi}$ denotes the Mellin transform of $\psi$.

Furthermore, for $X\ge 1$ and some constants $\beta_j$ depending on $\alpha_i$ we have
\begin{align}
\label{further} \Psi^{\pm}(X) = X \sum_{j=1}^{J} \beta_j \int_{0}^{\infty} \psi(x) e(\pm 3(xX)^{1/3} ) (xX)^{-j/3} dx + O(X^{(3-J)/3}).
\end{align}
\end{lemma}
\noindent Writing $s=\sigma+it$, by Stirling's approximation of the gamma function we have
\begin{align}
\label{Hbound} H^{\pm}(s)\ll_\sigma |t|^{3\sigma}.
\end{align}

We will need another result to estimate $\Psi^{\pm}(X)$ later.

\begin{lemma} \label{stationary} Suppose $\psi$ is a smooth, compactly supported function and $\gamma_1, \gamma_2> q^{-\epsilon}$ satisfy $ \gamma_1 \asymp \gamma_2$. Then we have
\begin{align}
\int_0^{\infty} \psi(x) e(2\gamma_1x^{1/2}\pm 3\gamma_2 x^{1/3}) dx \ll \gamma_1^{-1/2} (\| \psi \|_\infty +\| \psi' \|_\infty)q^{\epsilon}.
\end{align}
\end{lemma}
\proof
In the case of a plus sign, the required inequality follows by integrating by parts. In the case of a negative sign, the integrand has a stationary point at $x_0=(\gamma_2/\gamma_1)^6$. The required inequality follows by bounding the integral trivially in the range $|x-x_0|\le \gamma_1^{-1/2}$ and by integrating by parts in the range $|x-x_0|> \gamma_1^{-1/2}$. See section 9 of \cite{you} for a more precise treatment.
\endproof

\subsection{Other prerequisites}

We will need the following large sieve inequality. This may be found in Theorem 7.7 of \cite{iwakow}.
\begin{lemma}\label{lsieve}
Suppose that $\xi_m$ are some real numbers satisfying
\begin{align}
\| \xi_{m_1} - \xi_{m_2}\| \ge \delta,
\end{align}
for $m_1\neq m_2$, where $0<\delta\le 1/2$ and $\|x\|$ denotes the distance of $x$ from the closest integer. Then for any complex numbers $\lambda_n$ we have
\begin{align}
\sum_m \Big| \sum_{N\le n \le 2N} \lambda_n e(\xi_m n) \Big|^2 \le (\delta^{-1} + N ) \sum_{N\le n \le 2N} |\lambda_n|^2.
\end{align}
\end{lemma}

\section{Arithmetic side}

The goal of this section is to show that
\begin{align}
\frac{12}{q(k-1)} \sumh_{f\in H_k^{\star}(q)}L(\half,f\times g)L(\half, f) = \text{diagonal} + \text{off-diagonal} + O(q^{-1+\theta_1+\epsilon}),
\end{align}
where the diagonal equals the sum in the first line of (\ref{d-o5}) plus the sum in the first line of (\ref{d-o5t}) and the off-diagonal is the sum of the second and third lines of both (\ref{d-o5}) and (\ref{d-o5t}). The diagonal arises from the contribution of the $\delta(n,m)$ term on the arithmetic side of the trace formula. The rest of the arithmetic side gives rise to the off-diagonal and the error. In section \ref{sectiondiag} we show that the diagonal yields the main term of Theorem \ref{main}. In sections \ref{sectionoff1} and \ref{sectionoff2} we show that the off-diagonal falls into the error term of Theorem \ref{main}.

By the approximate functional equations and (\ref{hmult}), we have
\begin{align}
\frac{12}{q(k-1)}\sumh_{f\in H_k^{\star}(q)}L(\half,f\times g)L(\half, f) = \mathcal{M} + \widetilde{\mathcal{M}},
\end{align}
where
\begin{align}
\label{d-o} \mathcal{M} = &\sum_{n,m,r\ge 1}\frac{A(r,n)}{r\sqrt{nm}} V_1\Big(\frac{r^2n}{q^{3/2}}\Big)V_2\Big(\frac{m}{q^{1/2}}\Big)\Delta^{\star}_{k,q}(n,m)\\
&\nonumber  - i^{k}\sqrt{q}\sum_{n,m,r\ge 1}\frac{A(r,n)}{r\sqrt{nm}} V_1\Big(\frac{r^2n}{q^{3/2}}\Big)V_2\Big(\frac{m}{q^{1/2}}\Big)\Delta^{\star}_{k,q}(nq,m),
\end{align}
and
\begin{align}
\widetilde{\mathcal{M}} = &\sum_{n,m,r\ge 1}\frac{A(n,r)}{r\sqrt{nm}} \widetilde{V}_1\Big(\frac{r^2n}{q^{3/2}}\Big)V_2\Big(\frac{m}{q^{1/2}}\Big)\Delta^{\star}_{k,q}(n,m)\\
&\nonumber  - i^{k}\sqrt{q}\sum_{n,m,r\ge 1}\frac{A(n,r)}{r\sqrt{nm}} \widetilde{V}_1\Big(\frac{r^2n}{q^{3/2}}\Big)V_2\Big(\frac{m}{q^{1/2}}\Big)\Delta^{\star}_{k,q}(nq,m).
\end{align}
 By the trace formula, the first line of (\ref{d-o}) equals
\begin{align}
\label{d-o3}
&\sum_{n,r\ge 1}\frac{A(r,n)}{rn} V_1\Big(\frac{r^2n}{q^{3/2}}\Big)V_2\Big(\frac{n}{q^{1/2}}\Big)\\
&\nonumber+2\pi i^k  \sum_{n,m,r,c\ge 1} \frac{A(r,n)}{r\sqrt{nm}}\frac{S(n,m;cq)}{cq}J_{k-1}\Big(\frac{4\pi\sqrt{nm}}{cq}\Big)V_1\Big(\frac{r^2n}{q^{3/2}}\Big)V_2\Big(\frac{n}{q^{1/2}}\Big)\\
&\nonumber- \frac{12}{q(k-1)}\sumh_{f\in H^{\star}_k(1)} \sum_{\substack{n,m,r\ge 1}}  \frac{A(r,n)a_f(n)a_f(m)}{r\sqrt{nm}}V_1\Big(\frac{r^2n}{q^{3/2}}\Big)V_2\Big(\frac{m}{q^{1/2}}\Big)+O(q^{-1+\epsilon}).
\end{align}
The last line above is $\ll q^{-1+\epsilon}$ by Lemma \ref{vnote}. 

By separating the cases $q\nmid n$ and $q|n$, using (\ref{hmult}) and the identity $a_f(q)^2=q^{-1}$ for $f\in H_k^{\star}(q)$, we have that the second line of (\ref{d-o}) equals
\begin{align}
\label{d-o2}
&- i^{k}\sqrt{q}\sum_{\substack{n,m,r\ge 1\\(n,q)=1}}\frac{A(r,n)}{r\sqrt{nm}} V_1\Big(\frac{r^2n}{q^{3/2}}\Big)V_2\Big(\frac{m}{q^{1/2}}\Big)\Delta^{\star}_{k,q}(nq,m)\\
&\nonumber - \frac{i^k}{q} \sum_{n,m,r\ge 1}\frac{A(r,nq)}{r\sqrt{nm}} V_1\Big(\frac{r^2n}{q^{1/2}}\Big)V_2\Big(\frac{m}{q^{1/2}}\Big)\Delta^{\star}_{k,q}(n,m).
\end{align}
The last line of (\ref{d-o2}) is $\ll q^{-1+\theta_1+ \epsilon}$ by (\ref{vbound}) and (\ref{tracenote}).  By the trace formula, the first line of (\ref{d-o2}) equals
\begin{multline}
\label{d-o4} -2\pi \sqrt{q} \sum_{\substack{n,m,r,c\ge 1\\ (n,q)=1}} \frac{A(r,n)}{r\sqrt{nm}}\frac{S(nq,m,cq)}{cq}J_{k-1}\Big(\frac{4\pi\sqrt{nqm}}{cq}\Big)V_1\Big(\frac{r^2n}{q^{3/2}}\Big)V_2\Big(\frac{m}{q^{1/2}}\Big)\\
+O\Big(\sumh_{f\in H_k^{\star}(1)} \frac{a_f(q)}{q^{3/2}}\sum_{n,m,r\ge 1}   \frac{A(r,n)a_f(n)a_f(m)}{r\sqrt{nm}}V_1\Big(\frac{r^2n}{q^{3/2}}\Big)V_2\Big(\frac{m}{q^{1/2}} \Big) + q^{-3/2+\epsilon} \Big).
\end{multline}
The error term is $O(q^{-3/2+\epsilon})$ by Lemma \ref{vnote}. By (\ref{jbound}) and (\ref{vbound}), the contribution to (\ref{d-o4}) by the terms with $c>q^{1/2+\epsilon}$ or $m>q^{1/2+\epsilon}$ is $\ll q^{-100}$ for large enough $k$. Thus we may assume that $(c,q)=(m,q)=1$, so that we have the identity $S(nq,m,cq)= S(nq\overline{c},m\overline{c},q)S(n,m\overline{q};c)=-S(n\overline{q},m,c)$. We may also assume that $(rn,q)=1$ by (\ref{vbound}), so that $A(r,nq)=A(r,n)A(1,q)$. To the sum in (\ref{d-o4}) we may add the terms with $q|n$, incurring an error of
\begin{align}
2\pi\sqrt{q}& \sum_{n,m,r,c\ge 1}  \frac{A(r,nq)}{r\sqrt{nqm}}\frac{S(n,m;c)}{cq}J_{k-1}\Big(\frac{4\pi\sqrt{nq^{2}m}}{cq}\Big)V_1\Big(\frac{r^2n}{q^{1/2}}\Big)V_2\Big(\frac{m}{q^{1/2}}\Big)\\
\nonumber &\ll q^{-1+\theta_1+\epsilon} +  \sumh_{f\in H_k^{\star}(1)}  \frac{A(1,q)}{q} \sum_{n,m,r\ge 1}  \frac{A(r,n)a_f(n)a_f(m)}{r\sqrt{nm}}  V_1\Big(\frac{r^2n}{q^{1/2}}\Big)V_2\Big(\frac{m}{q^{1/2}}\Big)\\
\nonumber &\ll q^{-1+\theta_1+\epsilon},
\end{align}
where we used (\ref{ptf}) to evaluate the $c$-sum exactly and then used Lemma \ref{vnote}.

Gathering everything together, we have shown that
\begin{align}
\label{d-o5} \mathcal{M}=& \sum_{n,r\ge 1}\frac{A(r,n)}{rn} V_1\Big(\frac{r^2n}{q^{3/2}}\Big)V_2\Big(\frac{n}{q^{1/2}}\Big)\\
&\nonumber+2\pi i^k  \sum_{n,m,r,c\ge 1} \frac{A(r,n)}{r\sqrt{nm}}\frac{S(n,m;cq)}{cq}J_{k-1}\Big(\frac{4\pi\sqrt{nm}}{cq}\Big)V_1\Big(\frac{r^2n}{q^{3/2}}\Big)V_2\Big(\frac{m}{q^{1/2}}\Big)\\
&\nonumber+2\pi \sqrt{q} \sum_{n,m,r,c\ge 1} \frac{A(r,n)}{r\sqrt{nm}}\frac{S(n\overline{q},m,c)}{cq}J_{k-1}\Big(\frac{4\pi\sqrt{nqm}}{cq}\Big)V_1\Big(\frac{r^2n}{q^{3/2}}\Big)V_2\Big(\frac{m}{q^{1/2}}\Big)\\
&\nonumber+O(q^{-1+\theta_1+\epsilon}).
\end{align}
Similarly,
\begin{align}
\label{d-o5t} \widetilde{\mathcal{M}}=& \sum_{n,r\ge 1}\frac{A(n,r)}{rn} \widetilde{V}_1\Big(\frac{r^2n}{q^{3/2}}\Big)V_2\Big(\frac{n}{q^{1/2}}\Big)\\
&\nonumber+2\pi i^k  \sum_{n,m,r,c\ge 1} \frac{A(n,r)}{r\sqrt{nm}}\frac{S(n,m;cq)}{cq}J_{k-1}\Big(\frac{4\pi\sqrt{nm}}{cq}\Big)\widetilde{V}_1\Big(\frac{r^2n}{q^{3/2}}\Big)V_2\Big(\frac{m}{q^{1/2}}\Big)\\
&\nonumber+2\pi \sqrt{q} \sum_{n,m,r,c\ge 1} \frac{A(n,r)}{r\sqrt{nm}}\frac{S(n\overline{q},m,c)}{cq}J_{k-1}\Big(\frac{4\pi\sqrt{nqm}}{cq}\Big)\widetilde{V}_1\Big(\frac{r^2n}{q^{3/2}}\Big)V_2\Big(\frac{m}{q^{1/2}}\Big)\\
&\nonumber+O(q^{-1+\theta_1+\epsilon}).
\end{align}

\section{Diagonal}\label{sectiondiag}
In this section we evaluate the first lines of (\ref{d-o5}) and (\ref{d-o5t}). Together they form the diagonal contribution, which yields the main term of Theorem \ref{main}. 

\begin{lemma}
\begin{multline}
\label{diago} \sum_{n,r\ge 1}\frac{A(r,n)}{rn} V_1\Big(\frac{r^2n}{q^{3/2}}\Big)V_2\Big(\frac{n}{q^{1/2}}\Big) +\sum_{n,r\ge 1}\frac{A(n,r)}{rn} \widetilde{V}_1\Big(\frac{r^2n}{q^{3/2}}\Big)V_2\Big(\frac{n}{q^{1/2}}\Big)\\
= \frac{L(1,\tilde{g})L(1,g)}{\zeta(2)}\Big(1+\frac{\widetilde{G}_1(\half)}{G_1(\half)}\Big) + O(q^{-1}).
\end{multline}
\end{lemma}
\proof
By definition we have that the left hand side of (\ref{diago}) equals
\begin{multline}
\Big(\frac{1}{2\pi i}\Big)^2 \int_{(2)}\int_{(2)} \sum_{n,r\ge 1}\frac{A(r,n)}{r^{1+2s_1}n^{1+s_1+s_2}} q^{\frac{3}{2}s_1+\frac{1}{2}s_2}\frac{G_1(\half+s_1)G_2(\half+s_2)}{G_1(\half)G_2(\half)}\frac{ds_1ds_2}{s_1s_2}\\
+\Big(\frac{1}{2\pi i}\Big)^2 \int_{(2)}\int_{(2)} \sum_{n,r\ge 1}\frac{A(n,r)}{r^{1+2s_1}n^{1+s_1+s_2}} q^{\frac{3}{2}s_1+\frac{1}{2}s_2}\frac{\widetilde{G}_1(\half+s_1)G_2(\half+s_2)}{G_1(\half)G_2(\half)}\frac{ds_1ds_2}{s_1s_2}.
\end{multline}
We use Bump's identity (cf. Proposition 6.6.3 of \cite{gold})
\begin{align}
\sum_{n,r\ge 1}\frac{A(r,n)}{r^{s_1}n^{s_2}} = \frac{L(s_1,\tilde{g})L(s_2,g)}{\zeta(s_1+s_2)}
\end{align}
for $\Re(s_1),\Re(s_2)>1$, and then move the lines of integration to $\Re(s_1)=\Re(s_2)=-\half$. We pick up residues at $\Re(s_1)=0$ and $\Re(s_2)=0$, which form the main term. The integral on $\Re(s_1)=\Re(s_2)=-\half$ falls into the error term. Here we used the fact that $G_i(s)$ and $\widetilde{G}_i(s)$ decay exponentially on vertical lines, by Stirling's estimates.
\endproof

\section{Off-diagonal- Part 1}\label{sectionoff1}

For the second line of (\ref{d-o5}), we show that
\begin{multline}
\label{third} \sum_{n,m,r,c\ge 1} \frac{A(r,n)}{r\sqrt{nm}}\frac{S(n,m;cq)}{cq}J_{k-1}\Big(\frac{4\pi\sqrt{nm}}{cq}\Big)V_1\Big(\frac{r^2n}{q^{3/2}}\Big)V_2\Big(\frac{m}{q^{1/2}}\Big)\\
\ll q^{-1/4+\theta_2/2+\epsilon} + q^{-1/2+\theta_1+\epsilon}.
\end{multline}
A similar argument gives the same bound for the second line of  (\ref{d-o5}).
Observe that we may restrict our attention to the terms satisfying $c, r \asymp 1$, $n \asymp q^{3/2}$ and $m\asymp q^{1/2}$. The contribution by the terms not satisfying these conditions is $\ll q^{-100}$ for large enough $k$. This follows by (\ref{vbound}) and (\ref{jbound}). Thus (\ref{third}) follows from

\begin{lemma} For $r,c \asymp 1$, we have
\begin{align}
\label{off2}\sum_{m\asymp q^{1/2}}  \frac{1}{\sqrt{m}} V_2\Big(\frac{m}{q^{1/2}}\Big)  \sum_{ n\ge 1}A(r,n) S(n,m;cq) W_1\Big(\frac{n}{q^{3/2}}\Big)\ll q^{3/2+\theta_2/2+\epsilon}+q^{5/4+\theta_1+\epsilon},
\end{align}
where $W_1(x)=x^{-1/2}J_{k-1}(4\pi c^{-1} \sqrt{x m q^{-1/2}}) V_1(r^2 x) \vartheta_1(x)$ and $\vartheta_1$ is a smooth function compactly supported on $[q^{-\epsilon},q^{\epsilon}]$ and satisfying $\vartheta_1^{(B)}(x) \ll_B (q^{\epsilon})^B$ for integers $B\ge 0$.
\end{lemma}
\proof 

{\it First application of Voronoi:} We apply the Voronoi formula to the $n$-sum in (\ref{off2}) after writing $S(n,m;cq)=\summ_{h \bmod cq} e((n\overline{h}+mh)/cq)$. It suffices to show, since all other terms arising from the Voronoi formula are similar, that
\begin{multline}
\label{off3}\sum_{m\asymp q^{1/2}}  \frac{V_2(\frac{m}{q^{1/2}})}{\sqrt{m}}   cq \summ_{h \bmod cq} e(mh/cq) \sum_{\substack{n\ge1\\ l | cqr}} \frac{ A(n, l)}{nl} S(rh, n; qcr/l) W_2\Big(\frac{n}{q^{3/2}}\Big)\\
\ll q^{3/2+\theta_2/2+\epsilon}+q^{5/4+\theta_1+\epsilon},
\end{multline}
where
\begin{align}
\label{defw2} W_2(X)= X \int_{(\sigma)} (\pi^3X)^{-s} (l^2/c^3r)^{-s+1} H^+(s) \widetilde{W_1}(1-s) ds
\end{align}
for $\sigma>\theta_2$.
Note that $W_1^{(B)}(x)\ll_B (q^{\epsilon})^B$ so that by integrating by parts $B$ times we have for $s=\sigma+it$,
\begin{align}
 \widetilde{W_1}(1-s) = \int_{q^{-\epsilon}}^{q^\epsilon} x^{-s} W_1(x) dx \ll_B |t|^{-B} (q^{\epsilon})^B.
\end{align}
We can use this bound together with (\ref{Hbound}) to estimate $W_2(X)$. If $X>q^{\epsilon}$, we take $\sigma=\frac{1000}{\epsilon}$ in (\ref{defw2}) to see that $W_2(X)\ll X^{-2}q^{-100}$. If $ X \le q^{\epsilon}$, we take $\sigma=1$ to see that $W_2(X)\ll q^\epsilon$.
So the $n$-sum in (\ref{off3}) is essentially supported on $n<q^{3/2+\epsilon}$. We open the Kloosterman sum: $S(rh, n; qcr/l) = \summ_{u \bmod qcr/l} e((rhu+n\overline{u})l/qcr)$.

The contribution to (\ref{off3}) by the terms with $q|l$ is
\begin{align}
&\sum_{m\asymp q^{1/2}}  \frac{V_2(\frac{m}{q^{1/2}})}{\sqrt{m}}  c \summ_{h \bmod cq} e(mh/cq) \sum_{\substack{n\ge1\\ l | cr}} \frac{ A(n, l)}{nl} S(rh, n; cr/l) W_2\Big(\frac{n}{q^{3/2}}\Big)\\
\nonumber &=\sum_{m\asymp q^{1/2}}  \frac{V_2(\frac{m}{q^{1/2}})}{\sqrt{m}}    c \sum_{\substack{n\ge1\\ l | cr}} \frac{ A(n, ql)}{nl} W_2\Big(\frac{n}{q^{3/2}}\Big) \summ_{u \bmod cr/l} e\Big(\frac{nl\overline{u}}{cr}\Big)  \summ_{h \bmod cq} e\Big(\frac{h(m+uql)}{cq}\Big)\\
\nonumber &\ll q^{1/4+\theta_1+\epsilon},
\end{align}
since the innermost sum of the second line, a Ramanujan sum, is $\ll \sum_{d|(cq,m+qul)} d \ll q^{\epsilon}$ as $m\asymp q^{1/2}$.

Henceforth fix $l|cr$, so that $l\asymp 1$. Note that now we can quote a superior bound for $W_2(X)$ when $X\le q^{\epsilon}$ by taking $\sigma=\theta_2+\epsilon$ in (\ref{defw2}). We have
\begin{align}
\label{w2bound} W_2(X)\ll q^{\epsilon}X^{1-\theta_2-\epsilon}.
\end{align}
Define
\begin{align}
W_3(x)= x^{-1}W_2(x).
\end{align}
Exchanging the order of summation in (\ref{off3}), it is enough to show that
\begin{multline}
\label{off4}\sum_{m\asymp q^{1/2}}  \frac{V_2(\frac{m}{q^{1/2}})}{\sqrt{m}}  \summ_{u \bmod qcr/l}  \sum_{n\ge 1} A(n,l) e(n\overline{u}l/qcr) W_3 \Big( \frac{n}{q^{3/2}}\Big) \summ_{h \bmod cq} e(h(ul+m)/cq)\\
\ll q^{2+\theta_2/2+\epsilon}+q^{7/4+\theta_1+\epsilon}.
\end{multline}
The innermost sum above, a Ramanujan sum, equals
\begin{align}
\Big(\summ_{h \bmod q} e(h(ul+m)/q)\Big)\Big(\summ_{h \bmod c} e(h(ul+m)/c)\Big),
\end{align}
since $(c,q)=1$. 
Note that $\summ_{h \bmod q} e(h(ul+m)/q)$ equals $-1$ or $q-1$ according as $ul \not\equiv -m\bmod q$ or $ul \equiv -m \bmod q$ respectively. So the left hand side of (\ref{off4}) equals
\begin{align}
& -\sum_{m\asymp q^{1/2}}  \frac{V_2(\frac{m}{q^{1/2}})}{\sqrt{m}}  \sum_{n\ge 1} A(n,l) W_3\Big( \frac{n}{q^{3/2}}\Big)  \summ_{h \bmod c}  e(hm/c) S(n,qhr;qcr/l) \label{off5} \\
& +\sum_{m\asymp q^{1/2}}  \frac{V_2(\frac{m}{q^{1/2}})}{\sqrt{m}}  q\sum_{n\ge 1} A(n,l) W_3\Big( \frac{n}{q^{3/2}}\Big)  \summ_{h \bmod c}  e(hm/c) \summ_{\substack{u \bmod qcr/l\\ u\equiv -m \overline{l}\bmod q}}e\Big(\frac{uqhr+\overline{u}n}{qcr/l} \Big).\nonumber
\end{align}
Since $(cr/l,q)=1$, we have $S(n,qhr;qcr/l)=S(n\overline{q},hr;cr/l)S(0,n;q)$. This product of a Kloosterman sum and a Ramanujan sum is $\ll q^{1+\epsilon}$ if $q|n$ and $\ll q^{\epsilon}$ otherwise. In any case, the first line of (\ref{off5}) is $\ll q^{7/4+\theta_1+\epsilon}$. Now consider the second line. By the Chinese Remainder Theorem, for $(u,qcr/l)=1$ and $u\equiv -m\overline{l} \bmod q$, we can write $u=-m\overline{cr}(cr/l) + vq$, where $cr\overline{cr}\equiv 1 \bmod q$ and $(v,cr/l)=1$.  We have 
\begin{align}
&e\Big(\frac{uqhr}{qcr/l}\Big)=e\Big(\frac{vhrq}{cr/l}\Big),\\
&e\Big(\frac{n\overline{u}}{qcr/l}\Big)=e\Big(\frac{n\overline{ucr/l}}{q}\Big)e\Big(\frac{n\overline{uq}}{cr/l}\Big)=e\Big(\frac{-nl^2\overline{mcr}}{q}\Big)e\Big(\frac{n\overline{vq^2}}{cr/l}\Big).\nonumber
\end{align}
The proof of the lemma is now reduced to showing
\begin{multline}
\label{off06} \sum_{m\asymp q^{1/2}}  \frac{V_2(\frac{m}{q^{1/2}})}{\sqrt{m}} \sum_{n\ge 1} A(n,l) e\Big(\frac{-nl^2\overline{mcr}}{q}\Big) W_3 \Big( \frac{n}{q^{3/2}}\Big) \summ_{h \bmod c} e\Big(\frac{hm}{c}\Big)S(n\overline{q^2},hrq;cr/l)\\
\ll q^{1+\theta_2/2+\epsilon}.
\end{multline}
Our argument now will proceed differently according to the size of $n$ in (\ref{off06}). Thus we take a smooth partition of unity of $\mathbb{R}^{+}$ subordinate to a covering by dyadic intervals. Let $\vartheta_2(x)$ be a smooth function, which is compactly supported on $[1,2]$ and satisfies $\vartheta_2^{(B)}(x)\ll_B 1$. We need to show
\begin{multline}
\label{off006} \sum_{m\asymp q^{1/2}}  \frac{V_2(\frac{m}{q^{1/2}})}{\sqrt{m}} \sum_{n\ge 1} A(n,l) e\Big(\frac{-nl^2\overline{mcr}}{q}\Big) W_3 \Big( \frac{n}{q^{3/2}}\Big)  \vartheta_2 \Big(\frac{n}{N}\Big)  \\
 \summ_{h \bmod c} e\Big(\frac{hm}{c}\Big)S(n\overline{q^2},hrq;cr/l) \ll q^{1+\theta_2/2+\epsilon},
\end{multline}
for $N<q^{3/2+\epsilon}$.

$\bullet$ {\it Case 1: $q < N <q^{3/2+\epsilon}$}.

{\it Application of reciprocity:}
By the Chinese Remainder Theorem we have 
\begin{align}
e\Big(\frac{-nl^2\overline{mcr}}{q}\Big)= e\Big(\frac{-nl^2}{mcrq}\Big)e\Big(\frac{nl^2\overline{q}}{mcr}\Big).
\end{align}
Define
\begin{align}
W_4(x)= e\Big(\frac{-xq^{1/2}l^2}{mcr}\Big) W_3(x),
\end{align}
and note that $\frac{q^{1/2}l^2}{mcr}\asymp 1$. To establish (\ref{off006}) in the present case, one opens the Kloosterman sum and observes that it is enough to show
\begin{align}
\label{off7} \sum_{n\ge 1} A(n,l) e(n \overline{b}/d) W_4 \Big( \frac{n}{q^{3/2}} \Big) \vartheta_2 \Big(\frac{n}{N}\Big) \ll q^{3/4+\theta_2/2+\epsilon},
\end{align}
where $ d<q^{1/2+\epsilon}$ and $(b,d)=1$. Let
\begin{align}
W_5(x) = W_4(xNq^{-3/2})\vartheta_2(x).
\end{align}

{\it Second Application of Voronoi:} We apply the Voronoi formula to the left hand side of (\ref{off7}). It suffices to show, since all other terms arising from the Voronoi formula are similar, that
\begin{align}
\label{off8} d \sum_{\substack{n\ge 1\\ \ell | d l}} \frac{A(\ell,n)}{n\ell} S(lb, n ; d l/ \ell) W_6\Big( \frac{n}{d^3N^{-1}l\ell^{-2}}  \Big) \ll  q^{3/4+\theta_2/2+\epsilon},
\end{align}
where
\begin{align}
\label{w6} W_6(X)= X \int_{(\sigma)} (\pi^3 X)^{-s}  H^+(s) \widetilde{W_5}(1-s) ds
\end{align}
for $\sigma> \theta_2$.
We need to estimate $W_6(X)$. To this end we first note that by (\ref{w2bound}), we have
\begin{align}
W_5^{(B)}(x)\ll_B \Big(\frac{q^{3/2}}{N}\Big)^{\theta_2} q^{\epsilon} \ll_B q^{\theta_2/2+\epsilon}.
\end{align}
Integrating by parts $B$ times, we have for $s=\sigma+it$ the bound
\begin{align}
 \widetilde{W_5}(1-s) = \int_{1}^{2} x^{-s} W_5(x) dx &\ll_B |t|^{-B} q^{\theta_2/2+\epsilon} .
 \end{align}
Now we can estimate $W_6(X)$. If $X>q^{\epsilon}$ then by taking $\sigma=\frac{1000}{\epsilon}$ in (\ref{w6}), we see that $W_6(X)\ll q^{-100}X^{-2}$. If $X < q^{\epsilon}$ then by taking $\sigma=1$ in (\ref{w6}), we see that $W_6(X)\ll q^{\theta_2/2+\epsilon}$.  The proof is now reduced to showing
\begin{align}
\sum_{\substack{n< q^{1/2+\epsilon}\\ \ell | dl}} \frac{|A(\ell,n)|}{n\ell}| S(lb, n ; d l/ \ell)| \ll q^{1/4+\epsilon}.
\end{align}
By (\ref{weil}), $| S(lb, n ; d l/ \ell)| \ll (lb,n,dl/\ell)^{1/2} (dl) ^{1/2+\epsilon}\ll q^{1/4+\epsilon}$, since $(b,d)=1$. Using this and (\ref{avg}) completes the proof in this case.

$\bullet$ {\it Case 2:  $N \le q $}.  In this case we apply the large sieve estimate contained in Lemma \ref{lsieve}. To set up for this application, we need to separate $n$ and $m$ in $W_3 \big( \frac{n}{q^{3/2}}\big)$. Let 
\begin{align}
W_3(X,x)=- \int_{(\theta_2+\epsilon)} (\pi^3 X)^{-s} (l^2/c^3r)^{-s+1} H^{+}(s) \frac{x^{-s+3}}{(-s+1)(-s+2)(-s+3)} ds.
\end{align}
Note that 
\begin{align}
W_3(X) = \int_{q^{-\epsilon}}^{q^{\epsilon}}   W_1^{(3)}(x)  W_3(X,x) dx,
\end{align}
and $W_3(X,x)\ll X^{-\theta_2- \epsilon}$ for $X<1$ and $x\asymp 1$. We open the Kloosterman sum:
\begin{align}
S(n\overline{q^2},hrq;cr/l) = \summ_{u\bmod cr/l} e\Big(\frac{nu\overline{q^2}+hrq\overline{u}}{cr/l}\Big).
\end{align}
We define for $x,u\asymp 1$,
\begin{align}
\lambda_m = \frac{1}{\sqrt{m}}V_2\Big(\frac{m}{q^{1/2}}\Big) W_1^{(3)}(x) e\Big(\frac{hm}{c}\Big),
\end{align}
and
\begin{align}
\lambda_n =A(n,l) W_3\Big(\frac{n}{q^{3/2}},x\Big)\vartheta_2\Big(\frac{n}{N}\Big)e\Big(\frac{n u \overline{q^2}}{cr/l}\Big).
\end{align}
Note that $\lambda_m\ll m^{-1/2+\epsilon}$ and $\lambda_n\ll |A(n,l)| (q^{3/2}/N)^{\theta_2+\epsilon}.$ The left hand side of (\ref{off006}) equals
\begin{align}
\summ_{h\bmod c} e\Big(\frac{hlq\overline{u}}{c}\Big) \summ_{u \bmod cr/l} \int_{q^{-\epsilon}}^{q^{\epsilon}}  \sum_{m\asymp q^{1/2}} \lambda_m \sum_{n\asymp  N} \lambda_n e(-nl^2\overline{mcr}/q) dx.
\end{align}
Thus it is enough to show
\begin{align}
\label{avg2} \sum_{m\asymp q^{1/2}} \lambda_m \sum_{n\asymp  N} \lambda_n e(-nl^2\overline{mcr}/q)\ll q^{1+\theta_2/2+\epsilon}.
\end{align}
By the Cauchy-Schwarz inequality, we have that the left hand side of (\ref{avg2}) is less than
\begin{align}
q^{\epsilon} \Big( \sum_{m\asymp q^{1/2}} \Big| \sum_{n\asymp  N} \lambda_n e(-nl^2\overline{mcr}/q) \Big|^2 \Big)^{1/2}.
\end{align}
We may now apply Lemma \ref{lsieve}, with $\delta= q^{-1}$, to get that
\begin{align}
\sum_{m\asymp q^{1/2}} \Big| \sum_{n\asymp  N} \lambda_n e(-nl^2\overline{mcr}/q) \Big|^2  &\ll q  \sum_{n\asymp N} |A(n,l)|^2 \Big(\frac{q^{3/2}}{N}\Big)^{2\theta_2+\epsilon} \ll q^{2+\theta_2+\epsilon},
\end{align}
for $N\le q$.
\endproof

\section{Off-diagonal- Part 2}\label{sectionoff2}

For the third line of (\ref{d-o5}), we show that
\begin{align}
\label{claimthird} \sqrt{q} \sum_{n,m,r,c\ge 1} \frac{A(r,n)}{r\sqrt{nm}}\frac{S(n\overline{q},m,c)}{cq}J_{k-1}\Big(\frac{4\pi\sqrt{nqm}}{cq}\Big)V_1\Big(\frac{r^2n}{q^{3/2}}\Big)V_2\Big(\frac{m}{q^{1/2}}\Big)\ll q^{-1/8+\epsilon}.
\end{align}
A similar argument gives the same bound for the third line of  (\ref{d-o5t}).
By (\ref{jbound}) and (\ref{vbound}), we assume that $k$ is large enough so that we can restrict the $c$-sum to $c<r^{-1}q^{1/2+\epsilon}$. Taking a smooth partition of unity of $\mathbb{R}^+$ subordinate to a covering by dyadic intervals, let $\omega_1(x)$ and $\omega_2(x)$ be smooth functions, compactly supported on $[1,2]$ and satisfying $\omega_i^{(B)}(x)\ll_B 1$. To establish (\ref{claimthird}), it suffices to prove the following:

\begin{lemma} Let $cr<q^{1/2+\epsilon}, N<q^{3/2+\epsilon}$ and $ M<q^{1/2+\epsilon}$. We have
\begin{align}
\label{r1}  \sum_{n,m\ge 1} A(r,n)\frac{S(n\overline{q},m;c)}{c} \Omega_1\Big(\frac{n}{N},\frac{m}{M}\Big) \ll  \frac{(NM)^{1/2}q^{3/8+\epsilon}}{cr},
\end{align}
where
\begin{align}
\Omega_1(x,y) = \omega_1(x) \omega_2(y) J_{k-1}\Big( \frac{4\pi\sqrt{xyNMq^{-1}}}{c} \Big)  V_1\Big(\frac{xr^2N}{q^{3/2}}\Big) V_2\Big(\frac{yM}{q^{1/2}}\Big).
\end{align}
\end{lemma}
\proof
We write $S(n\overline{q},m,c) = \summ_{h \bmod c} e(n\overline{qh}/c)e(mh,c)$ and apply the Voronoi formula to the $n$-sum. It suffices to show that
\begin{align}
\label{r2} \sum_{\substack{n,m\ge 1\\ l|rc}} \frac{A(n,l)}{nl} \Psi_{\frac{m}{M}}^{\pm} \Big(\frac{nNl^2}{c^3r}\Big) \summ_{h \bmod c}   e(hm/c)S(rhq,\pm n,rc/l) \ll \frac{(NM)^{1/2}q^{3/8+\epsilon}}{cr},
\end{align}
where
\begin{align}
\label{rpsi} \Psi_y^{\pm} (X)= \int_{(\sigma)} X^{1-s} H^{\pm}(s)\widetilde{\Omega}_1(1-s)  ds, 
\end{align}
for $\sigma>\theta_2$ and
\begin{align}
\widetilde{\Omega}_1(1-s)= \int_1^2 x^{-s} \Omega_1(x,y) dx.
\end{align}
By (\ref{jbound2}) we have
\begin{align}
 \frac{\partial^{B}}{\partial x^{B}} \Omega_1(x,y) \ll_B  \Big( \frac{\sqrt{NMq^{-1}}}{c} \Big)^{B-1/2}(q^\epsilon)^B.
\end{align}
Thus, writing $s=\sigma+it$, we have by integration by parts $B$ times,
\begin{align}
\label{bou} \widetilde{\Omega}_1(1-s) \ll_B |t|^{-B}  \Big( \frac{\sqrt{NMq^{-1}}}{c} \Big)^{B-1/2}(q^\epsilon)^B.
\end{align}
Using this bound with $B=\lfloor 3\sigma +10  \rfloor$ and (\ref{Hbound}), we have
\begin{align}
\Psi_y^{\pm}\Big(\frac{nNl^2}{c^3r}\Big) \ll_\sigma q^{100} \Big( \frac{nl^2N}{c^3r} \Big)^{-\sigma}  \Big( \frac{\sqrt{NMq^{-1}}}{c} \Big)^{3\sigma} \ll_\sigma q^{100} (nl^2)^{-\sigma} \Big( \frac{r^2NM^3}{q^{3}} \Big)^{\sigma/2}. 
\end{align}
We see by taking $\sigma$ large enough that if $r^2N<q^{3/2-\epsilon}$ or $M<q^{1/2-\epsilon}$ then the lemma is easily proved. Henceforth assume that $r^2N \asymp q^{3/2}$ and $M \asymp q^{1/2}$. We also see by taking $\sigma$ large enough that we may restrict the sum in (\ref{r2}) to $n,l\asymp 1$, up to an error of $O(q^{-100})$. Thus it suffices to show that
\begin{align}
\label{r3} \sum_{m\ge 1}  \Psi_\frac{m}{M}^{\pm} \Big(\frac{nNl^2}{c^3r}\Big) \summ_{h \bmod c}   e(hm/c)S(rhq,\pm n,rc/l) \ll \frac{q^{11/8+\epsilon}}{cr},
\end{align}
for fixed $n, l\asymp 1$ with $l|cr$.
Opening the Kloosterman sum, the left hand side equals
\begin{align}
\label{r4} \sum_{m\ge 1}  \Psi_{\frac{m}{M}}^{\pm} \Big(\frac{nNl^2}{c^3r}\Big)  \summ_{u \bmod rc/l}e(\pm nl\overline{u}/rc)   \summ_{h \bmod c} e(h(m+uql)/c).
\end{align}
The innermost sum is a Ramanujan sum and equals $\sum_{d|c}\mu(c/d)d\delta_{m\equiv -uql \bmod d}$, where the delta symbol equals $1$ if the congruence is satisfied and $0$ otherwise. Exchanging the order of summation above, it suffices to show
\begin{align}
\label{r5} \summ_{u \bmod rc/l}e(\pm nl\overline{u}/rc) \sum_{d|c}\mu(c/d) d \sum_{\substack{m\ge 1\\m\equiv -uql \bmod d}} \Psi_{\frac{m}{M}}^{\pm} \Big(\frac{nNl^2}{c^3r}\Big)  \ll \frac{q^{11/8+\epsilon}}{cr}.
\end{align}

We need to estimate $\Psi_y^{\pm}(X)$ for $X=\frac{nNl^2}{c^3r}\asymp (\frac{q^{1/2}}{cr})^3$. We claim that
\begin{align}
\label{claim1} \Psi_y^{\pm}(X) \ll  \frac{q^{1/2+\epsilon}}{cr} .
\end{align}
This can be seen as follows. If $cr\asymp q^{1/2}$ then we take $\sigma=1$ in (\ref{rpsi}) and use the bound (\ref{bou}) with $N\asymp q^{3/2}/r^2$ to prove the claim.
If $cr<q^{1/2-\epsilon}$ then $X>1$. We write $J_{k-1}(x) = \Re( e^{ix}\mathcal{J}(x))$ as in (\ref{jlong}), and let
\begin{align}
\Omega_2(x,y) = \omega_1(x) \omega_2(y) \mathcal{J} \Big( \frac{4\pi\sqrt{xyNMq^{-1}}}{c} \Big)  V_1\Big(\frac{xr^2N}{q^{3/2}}\Big) V_2\Big(\frac{yM}{q^{1/2}}\Big).
\end{align}
By (\ref{further}) we have
\begin{align}
 \Psi_y^{\pm}(X) = X \sum_{j=1}^{1000/\epsilon} \beta_j \int_{0}^{\infty} \frac{\Re( \Omega_2(x,y) e(\frac{2(xyNM)^{1/2}}{q^{1/2}c})) e(\pm 3(xX)^{1/3}  )}{(xX)^{j/3}} dx + O(q^{-10}).
\end{align}
Now using Lemma \ref{stationary} with $\gamma_1 = \frac{(xyNM)^{1/2}}{q^{1/2}c}\asymp \frac{q^{1/2}}{cr}$ and $\gamma_2 = X^{1/3}\asymp \frac{q^{1/2}}{cr}$, and the bound $\mathcal{J} \big( \frac{4\pi\sqrt{xyNMq^{-1}}}{c} \big) \ll \big(\frac{q^{1/2}}{cr}\big)^{-1/2}$, we get (\ref{claim1}).
Similarly we have for the derivates,
\begin{align}
\label{claim2}  \frac{\partial^{B}}{\partial y^{B}} \Psi_y^{\pm} (X) \ll_B \Big( \frac{q^{1/2+\epsilon}}{cr} \Big)^{B+1}
\end{align}
for $X\asymp (\frac{q^{1/2}}{cr})^3$.
The proof of (\ref{r5}) proceeds according to two cases.

$\bullet$ {\it Case 1: $ cr \le q^{3/8}$}.  

\noindent Bounding left hand side of (\ref{r5}) absolutely, we find that it is less than
\begin{align}
 \sum_{u \bmod rc/l} \sum_{d|c} d \frac{M}{d} \frac{q^{1/2+\epsilon}}{cr} \ll \frac{q^{11/8+\epsilon}}{cr}.
\end{align}

$\bullet$ {\it Case 2: $q^{3/8}< cr <q^{1/2+\epsilon}$}.  

\noindent We evaluate the innermost sum in (\ref{r5}) using Poisson summation:
\begin{align}
\sum_{\substack{m\ge 1\\m\equiv -uql \bmod d}}\Psi_{\frac{m}{M}}^{\pm} \Big(\frac{nNl^2}{c^3r}\Big) = \frac{M}{d} \sum_{m=-\infty}^{\infty} e(muql/d) \widehat{\Psi}_{\frac{m}{d/M}}^{\pm} \Big(\frac{nNl^2}{c^3r}\Big),
\end{align}
where 
\begin{align}
\widehat{\Psi}_{Y}^{\pm} ( X ) = \int_{-\infty}^{\infty} \Psi_{y}^{\pm} (X) e(yY) dy = \int_{1}^{2} \Psi_{y}^{\pm} (X) e(yY) dy.
\end{align}
So (\ref{r5}) reduces to proving
\begin{align}
\label{r6} \sum_{m=-\infty}^{\infty} \widehat{\Psi}_{\frac{m}{d/M}}^{\pm} \Big(\frac{nNl^2}{c^3r}\Big)   S(\pm n, mqrc/d ;cr/l) \ll \frac{q^{7/8+\epsilon}}{cr},
\end{align}
where $d|c$. By (\ref{claim1}) we have for $X=\frac{nNl^2}{c^3r}$,
\begin{align}
\label{hatbound} \widehat{\Psi}_{Y}^{\pm} ( X ) \ll \frac{q^{1/2+\epsilon}}{cr}.
\end{align}
Integrating by parts and using (\ref{claim2}) gives the bound
\begin{align}
\frac{\partial^{B}}{\partial Y^{B}} \widehat{\Psi}_Y^{\pm} (X)\ll (q^{\epsilon})^B  \Big( \frac{q^{1/2}}{cr} \Big)^{B+1} Y^{-B},
\end{align}
Thus the $m$-sum in (\ref{r6}) is essentially supported on $|m|<\frac{d}{M}\frac{q^{1/2+\epsilon}}{cr}<q^{1/8+\epsilon}$. Using Weil's bound for the Kloosterman sum and the bound (\ref{hatbound}), the left hand side of (\ref{r6}) is less than
\begin{align}
\sum_{|m|<q^{1/8+\epsilon}} \frac{q^{1/2}}{cr} \sum_{d|c} q^{1/4+\epsilon}\ll  \frac{q^{7/8+\epsilon}}{cr}.
\end{align}

\endproof

{\bf Acknowledgments.} While working on this paper, I was supported by a grant from the European Research Council (grant agreement number 258713). I am grateful to Prof. Valentin Blomer for inviting me to G\"{o}ttingen and for the many helpful discussions I have had with him regarding this paper. I am thankful to Prof. Matthew Young for working with me on related topics. I would also like to thank Dr. Sheng-Chi Liu for providing me with a preprint of his article \cite{liu}.

\bibliographystyle{amsplain}

\bibliography{simultaneous-nonvanishing-2}

\end{document}